\documentclass[12pt,a4paper]{article}
\usepackage{amsfonts,amssymb,amsmath,amsthm,cases}
\usepackage[colorlinks=true,linkcolor=blue,citecolor=red, breaklinks = true]{hyperref}
\newtheoremstyle{thmstyle}{}{}{}{}{\bf}{}{ }{}
\theoremstyle{thmstyle}
\numberwithin{equation}{section}

\usepackage{amssymb}
\usepackage{graphics}
\usepackage{graphicx}
\usepackage{hyperref}
\usepackage{listings}
\usepackage{mathtools}
\usepackage{maple}
\usepackage[utf8]{inputenc}  
\usepackage[svgnames]{xcolor}
\usepackage{amsmath}
\usepackage{breqn}
\usepackage{textcomp}
\begin{document}
\title{{Bicomplex Mittag-Leffler Distribution}}
\author{Dharmendra Kumar Singh$^{1}$ and Chinmay Sharma$^{2}$\\Department of Mathematics\\
	School of Basic Sciences, UIET\\CSJM University, Kanpur\\ Email: drdksinghabp@gmail.com}

\date{}
\maketitle
\begin{center}
\textbf{Abstract}
\end{center}
Probability distribution theory helps in studying the impact of various dimensions in life while the Mittag-Leffler function and bicomplex are used in electromagnetism, quantum mechanics, and signal theory. Considering
the importance of both, the purpose of this paper is to introduce bicomplex Mittag-Leffler distribution theory with the help of the bicomplex Mittag-Leffler function. Moreover, it also tells us about the moment-generating
function, the first four moments, the mean, and the variance of this endeavour.

\section{Introduction}
As the importance of the Mittag-Leffler function is increasing in physical applications, many researchers are studying various extensions and applications of the Mittag-Leffler function. Some of them are the Incomplete Mittag-Leffler Function (Singh and Porwal 2013), a generalization of the Mittag-Leffler function and its properties ( Shukla and Prajapati 2007), and  A singular integral equation with a generalized Mittag-Leffler function in the kernel (Prabhakar 1971).\\
The Mittag-Leffler function was first proposed by Swedish mathematician Gosta Mittag-Leffler in 1903 (Mittag-Leffler 1903).
\begin{equation*}
  {\mathbb{E}}_{\alpha}(z) = \sum_{k=0}^{\infty}\frac{z^{k}}{\Gamma(\alpha k+1)}, \indent z, \alpha \in \mathbb{C}
\end{equation*}
is called the Mittag Leffler function of order $\alpha$. The Mittag-Leffler function is the direct generalization of an exponential function
\begin{equation*}
    e^{z} = \sum_{k=0}^{\infty}\frac{z^{k}}{\Gamma(k+1)}.
\end{equation*}
The first generalization of the Mittag-Leffler function was given by A. Wiman ( Wiman 1905)
\begin{equation*}
    {\mathbb{E}}_{\alpha,\beta}(z) = \sum_{k=0}^{\infty}\frac{z^{k}}{\Gamma(\alpha k+\beta)}, \indent z, \alpha, \beta \in \mathbb{C}, Re(\alpha) > 0, Re(\beta) > 0
\end{equation*}
is a two-parameter Mittag-Leffler function.\\*
\textbf{Bicomplex Number}\\
In terms of real components, the set of bicomplex numbers is defined as
\begin{equation*}
  \mathbb{T} = \{ \xi : \xi = x_{0} + i_{1}x_{1} + i_{2}x_{2} + jx_{3} \;| \; x_{0}, x_{1}, x_{2}, x_{3} \in \mathbb{R} \},
  \end{equation*}
and in terms of complex components, it can be written as
\begin{equation*}
  \mathbb{T} = \{\xi : \xi = z_{1} + i_{2}z_{2} \; | \; z_{1}, z_{2} \in \mathbb{C} \}
\end{equation*}
\begin{equation*}
    x_{0} = Re(\xi), \; x_{1} = Im_{i_{1}}(\xi), \; x_{2} = Im_{i_{2}}(\xi), \; x_{3} = Im_{j}(\xi).
\end{equation*}
\textbf{Null Cone}\\ The set of all zero divisors is called a null cone defined as
\begin{align*}
   \mathbb{NC} = \{ z_{1} + i_{2}z_{2} \;|\; z_{1}^{2} + z_{2}^{2} = 0 \}.
\end{align*}
Two non trivial zero divisors in $ \mathbb{T}$, denoted by $e_{1}$ and $e_{2}$ are defined as follows
\begin{equation*}
    e_{1} = \frac{1 + i_{1}i_{2}}{2} = \frac{1 + j}{2},
\end{equation*}
\begin{equation*}
    e_{2} = \frac{1 - i_{1}i_{2}}{2} = \frac{1 - j}{2},
\end{equation*}
\begin{equation}
    e_{1}.e_{2} = 0,\; e_{1} + e_{2} = 1,\; e_{1}^{2} = e_{1},\; e_{2}^{2} = e_{2}.
\end{equation}
\textbf{Idempotent Representation}\\ Every element $\xi \in \mathbb{T}$ has unique idempotent representation in terms of $e_{1}$ and $e_{2}$ defined by
\begin{equation}
    \xi = z_{1} + i_{2}z_{2} = \xi_{1}e_{1} + \xi_{2}e_{2},
    \end{equation}
    where
    \begin{equation*}
        \xi_{1} = (z_{1} - i_{1}z_{2}),\; \xi_{2} = (z_{1} + i_{1}z_{2}).
    \end{equation*}
  \textbf{Bicomplex Mittag-Leffler Function}\\ The bicomplex one parameter Mittag-Leffler function is defined as
    \begin{equation}
  {\mathbb{E}}_{\alpha}(\xi) = \sum_{k=0}^{\infty}
        \frac{\xi^{k}}{\Gamma(\alpha k + 1)},\indent \xi, \alpha\in \mathbb{T}, |Im_{j}(\alpha)| < Re(\alpha).
    \end{equation}
    By using idempotent representation we can write $ E_{\alpha}(\xi)$ as follows
    \begin{equation}
       {\mathbb{E}}_{\alpha}(\xi) = \sum_{k=0}^{\infty}
        \frac{\xi_{1}^{k}}{\Gamma(\alpha_{1} k + 1)}e_{1} +  \sum_{k=0}^{\infty}
        \frac{\xi_{2}^{k}}{\Gamma(\alpha_{2} k + 1)}e_{2} =  E_{\alpha_{1}}(\xi_{1})e_{1} + E_{\alpha_{2}}(\xi_{2})e_{2}
    \end{equation}
    where
    \begin{equation}
     \xi =  \xi_{1}e_{1} + \xi_{2}e_{2},
     \end{equation}
\begin{equation}
          \alpha = \alpha_{1}e_{1} + \alpha_{2}e_{2}.
    \end{equation}
Researchers have also studied various distributions involving the Mittag-Leffler function. Mittag-Leffler Functions and Their Applications (Haubold et al. 2011), On Mittag Leffler Functions and Related Distributions (Pillai 1990), On Mittag-Leffler type Poisson Distribution (Porwal and Dixit 2017), An Overview of Generalized Gamma Mittag-Leffler Model and its Applications (Nair 2015), The role of the Mittag-Leffler function in fractional Modeling ( Rogosin 2015).\\
\indent Researchers have also explored Exotic Probability Theory, a branch of probability theory in which we study the probabilities that are outside the interval $ [0,1]$. Mostly the work has been done where the probabilities are either negative or complex numbers. As Exploring Exotic Probability Theory (Yannan Lin 2020-21), Quantum Mechanics as an Exotic Probability Theory (Saul Youssef 1995).\\
\indent This motivated us to define a probability density function that deals with bicomplex probabilities.
In this paper, we have introduced a new distribution named \textbf{Bicomplex Mittag-Leffler Distribution} that involves the Bicomplex Mittag-Leffler Function of one parameter (Agrwal et al. 2022). We have derived the Moment Generating Function, first four moments, Mean, and Variance of this distribution using the properties of bicomplex numbers (Agrwal et al. 2022).\\
    \textbf{Bicomplex Mittag-Leffler Distribution}

Consider the probability density function
\begin{equation}
    F(\xi) = (1+a)e^{-\xi}{\mathbb{E}}_{\alpha}(-a\xi^{\alpha}),
  \end{equation}
 where
        \begin{equation*}
  {\mathbb{E}}_{\alpha}(-a\xi^{\alpha}) = \sum_{k=0}^{\infty}
        \frac{(-a \xi^{\alpha})^{k}}{\Gamma(\alpha k + 1)},
    \end{equation*}
Here $\xi, \alpha\in \mathbb{T}$. $\xi=z_{1}+i_{2}z_{2}$ and $|Im_{j}(\alpha)|<\Re(\alpha).$\\
\textbf{Moment Generating Function}\\ If $\xi$ is a continuous random variable with probability density function
$F(\xi)$, then the moment generating function is defined as follows
\begin{equation}
    M_{\xi}(t) = E(e^{\xi t}).
\end{equation}
\begin{equation}
    M_{\xi}(t) = \int_{-\infty}^{\infty}e^{\xi t}F(\xi)d\xi.
\end{equation}
\textbf{rth moment}\\ If $\xi$ is a continuous random variable with probability density function
$F(\xi)$, then the rth moment about the origin  is defined as follows
\begin{equation}
    \mu_{r}^{'} = E(\xi^{r}) =  \int_{-\infty}^{\infty}\xi^{r}F(\xi)d\xi.
\end{equation}
\section{Results}
\subsection{Probability Density Function}
The function defined in (1.7) is a probability density function since

\begin{equation*}
    \int_{0}^{\infty} F(\xi)d\xi = 1.
\end{equation*}
\textbf{Proof}
\begin{equation*}
   \int_{0}^{\infty} F(\xi)d\xi = \int_{0}^{\infty}(1+a)e^{-\xi}{\mathbb{E}}_{\alpha}(-a\xi^{\alpha})d\xi
\end{equation*}
Using idempotent representation
  \begin{equation*}
    = (1+a)\left[\left(\int_{0}^{\infty}e^{-\xi_{1}}{\mathbb{E}}_{\alpha}(-a\xi_{1}^{\alpha})d\xi_{1}\right) e_{1} + \left(\int_{0}^{\infty}e^{-\xi_{2}}{\mathbb{E}}_{\alpha}(-a\xi_{2}^{\alpha})d\xi_{2}\right)e_{2} \right]
     \end{equation*}
 \begin{equation*}
     = \sum_{k=0}^{\infty}(1+a)\left[\left( \int_{0}^{\infty}e^{-\xi_{1}}\frac{(-a)^{k}\xi_{1}^{\alpha k} } {\Gamma(\alpha k+1)}d\xi_{1}\right)e_{1} + \left(\int_{0}^{\infty}e^{-\xi_{2}}\frac{(-a)^{k}\xi_{2}^{\alpha k} }{\Gamma(\alpha k+1)}d\xi_{2}\right)e_{2}\right]
           \end{equation*}
 \begin{equation*}
     = \sum_{k=0}^{\infty}(1+a)\left[\left( \frac{(-a)^{k}} {\Gamma(\alpha k+1)}\Gamma(\alpha k+1)\right)e_{1} + \left(\frac{(-a)^{k} }{\Gamma(\alpha k+1)}\Gamma(\alpha k+1) \right)e_{2}\right]
           \end{equation*}
           \begin{equation*}
              = (1+a)\sum_{k=0}^{\infty}(-a)^{k}(e_{1} + e_{2}).
           \end{equation*}
           Using(1.1)
           \begin{equation*}
               \int_{0}^{\infty} F(\xi)d\xi = 1.
           \end{equation*}
           
           We have relaxed the condition $ F(\xi) \ge 0    $ since we are dealing with the  bicomplex probabilities that lie outside the region $ [0,1]$.
\subsection{Moment Generating Function}
The moment generating function of the probability density function defined in (1.7) is given by
\begin{equation*}
    M_{\xi}(t) = E(e^{\xi t})
\end{equation*}
\begin{equation*}
    = \int_{0}^{\infty}e^{\xi t}F(\xi)d\xi
\end{equation*}
\begin{equation*}
 = \int_{0}^{\infty}e^{\xi t}(1+a)e^{-\xi}{\mathbb{E}}_{\alpha}(-a\xi^{\alpha})d\xi
\end{equation*}
\begin{equation*}
= (1+a)\int_{0}^{\infty}e^{-\xi(1-t)}{\mathbb{E}}_{\alpha}(-a\xi^{\alpha})d\xi.
    \end{equation*}
 Using idempotent representation
    \begin{equation*}
    M_{\xi}(t) = (1+a)\left[\left(\int_{0}^{\infty}e^{-\xi_{1}(1-t)}{\mathbb{E}}_{\alpha}(-a\xi_{1}^{\alpha})d\xi_{1}\right) e_{1} + \left(\int_{0}^{\infty}e^{-\xi_{2}(1-t)}{\mathbb{E}}_{\alpha}(-a\xi_{2}^{\alpha})d\xi_{2}\right)e_{2} \right]
     \end{equation*}
      \begin{equation*}
     = \sum_{k=0}^{\infty}(1+a)\left[\left( \int_{0}^{\infty}e^{-\xi_{1}(1-t)}\frac{(-a)^{k}\xi_{1}^{\alpha k} } {\Gamma(\alpha k+1)}d\xi_{1}\right)e_{1} + \left(\int_{0}^{\infty}e^{-\xi_{2}(1-t)}\frac{(-a)^{k}\xi_{2}^{\alpha k} }{\Gamma(\alpha k+1)}d\xi_{2}\right)e_{2}\right]
           \end{equation*}
           \begin{equation*}
              = \sum_{k=0}^{\infty}(1+a)(-a)^{k}\left[\left(\frac{1}{ \Gamma(\alpha k+1)} \frac{\Gamma(\alpha k+1)}{(1-t)^{\alpha k+1}}\right)e_{1} + \left(\frac{1}{ \Gamma(\alpha k+1)} \frac{\Gamma(\alpha k+1)}{(1-t)^{\alpha k+1}}\right)e_{2}\right]
           \end{equation*}
 \begin{equation*}
         = \sum_{k=0}^{\infty}(1+a)(-a)^{k}\left[\left(\frac{1}{(1-t)^{\alpha k+1}}\right)e_{1} + \left( \frac{1}{(1-t)^{\alpha k+1}}\right)e_{2}\right]
           \end{equation*}
           \begin{align*}
               =(1+a) \left[ \left(\frac{(1-t)^{\alpha -1}}{(1-t)^{\alpha}+a} \right) e_{1} + \left( \frac{(1-t)^{\alpha -1}}{(1-t)^{\alpha}+a}\right) e_{2} \right].
           \end{align*}
            ( where  $\left| \frac{a}{(1-t)^{\alpha }} \right| < 1 $)
           \begin{equation*}
 M_{\xi}(t)  = \frac{(1+a)(1-t)^{\alpha - 1}}{a + (1 - t)^{\alpha}}(e_{1} + e_{2}).
   \end{equation*}
 Using(1.1)  
   \begin{equation}
 M_{\xi}(t)  = \frac{(1+a)(1-t)^{\alpha - 1}}{a + (1 - t)^{\alpha}}.
   \end{equation}
\mapleplot{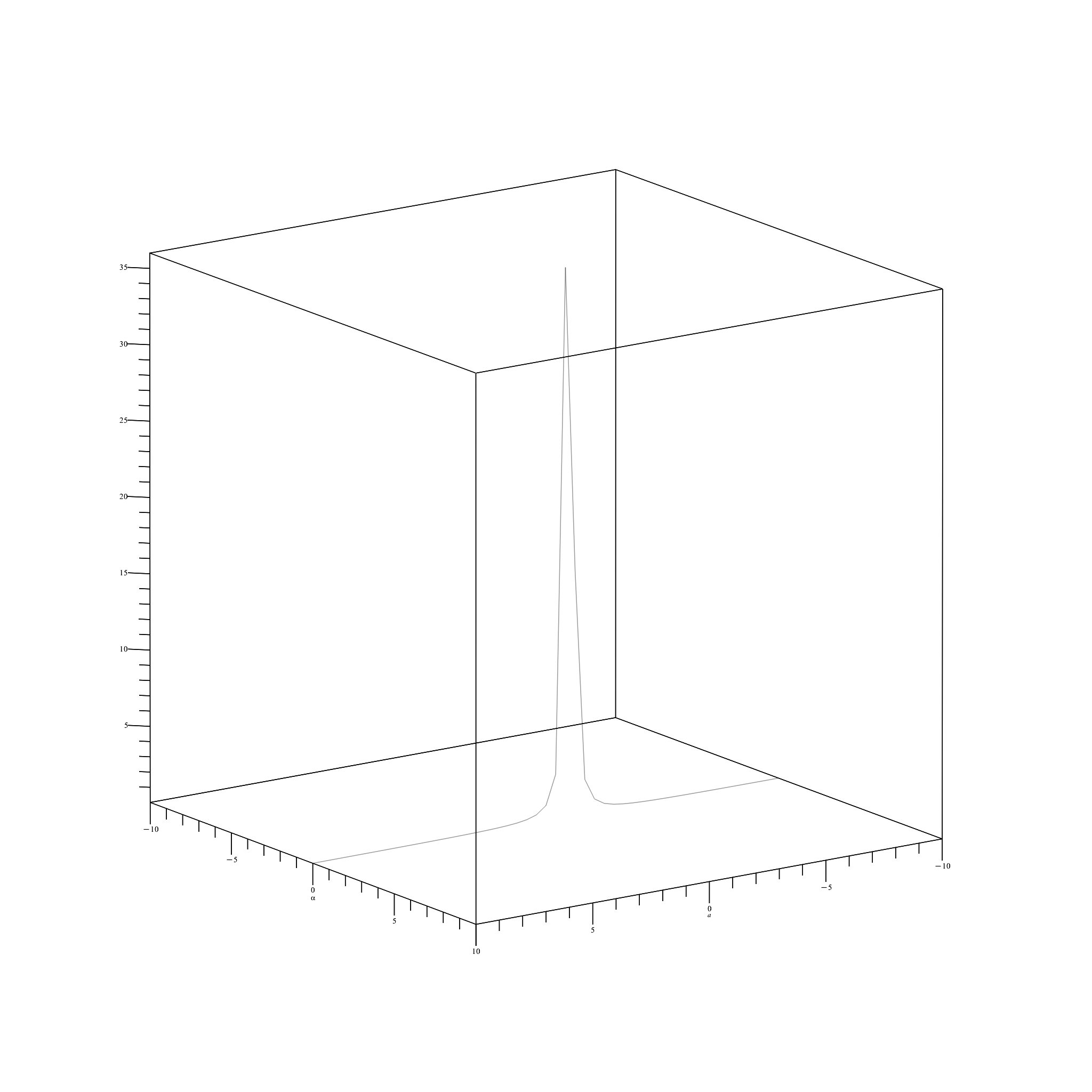}
\vspace{-2cm}
\begin{center}
	\footnotesize{Momont Generating Function}
\end{center}
 \subsection{First Moment About the Origin}
 \begin{equation*}
\mu_{1}^{'} = E(\xi)
    \end{equation*}
    \begin{equation*}
     =\int_{0}^{\infty} \xi F(\xi)d\xi
    \end{equation*}
\begin{equation*}
    =\int_{0}^{\infty}\xi(1+a)e^{-\xi}{\mathbb{E}}_{\alpha}(-a\xi^{\alpha})d\xi.
\end{equation*}
Using idempotent Representation
  \begin{equation*}
    \mu_{1}^{'} = (1+a)\left[\left(\int_{0}^{\infty}\xi_{1}e^{-\xi_{1}}{\mathbb{E}}_{\alpha}(-a\xi_{1}^{\alpha})d\xi_{1}\right) e_{1} + \left(\int_{0}^{\infty}\xi_{2}e^{-\xi_{2}}{\mathbb{E}}_{\alpha}(-a\xi_{2}^{\alpha})d\xi_{2}\right)e_{2} \right]
     \end{equation*}
     \begin{equation*}
  = \sum_{k=0}^{\infty}(1+a)\left[\left( \int_{0}^{\infty}\xi_{1}e^{-\xi_{1}}\frac{(-a)^{k}\xi_{1}^{\alpha k} } {\Gamma(\alpha k+1)}d\xi_{1} \right)e_{1} + \left(\int_{0}^{\infty}\xi_{2}e^{-\xi_{2}}\frac{(-a)^{k}\xi_{2}^{\alpha k} }{\Gamma(\alpha k+1)}d\xi_{2}\right)e_{2}\right]
           \end{equation*}
  \begin{equation*}
     = \sum_{k=0}^{\infty}(1+a)\left[\left(\frac{(-a)^{k}}{\Gamma(\alpha k+1)} \Gamma(\alpha k+2)\right)e_{1}+ \left(\frac{(-a)^{k}}{\Gamma(\alpha k+1)} \Gamma(\alpha k+2)\right)e_{2}\right]
  \end{equation*}
  \begin{equation*}
     = (1+a)\sum_{k=0}^{\infty}[{(-a)^{k}(\alpha k+1)} e_{1} + {(-a)^{k}(\alpha k+1)} e_{2}]
  \end{equation*}
  \begin{equation*}
     = (1+a)\left[\left(\frac{a+1-a\alpha}{(a+1)^{2}}\right)e_{1} + \left(\frac{a+1-a\alpha}{(a+1)^{2}}\right)e_{2}\right]
  \end{equation*}
  \begin{equation*}
      = \left(1- \frac{a\alpha}{a+1}\right)( e_{1} +  e_{2}).
  \end{equation*}
  
  Using (1.1) 
           \begin{equation}
             \mu_{1}^{'}  = 1- \frac{a\alpha}{a+1}.
           \end{equation}
       \mapleplot{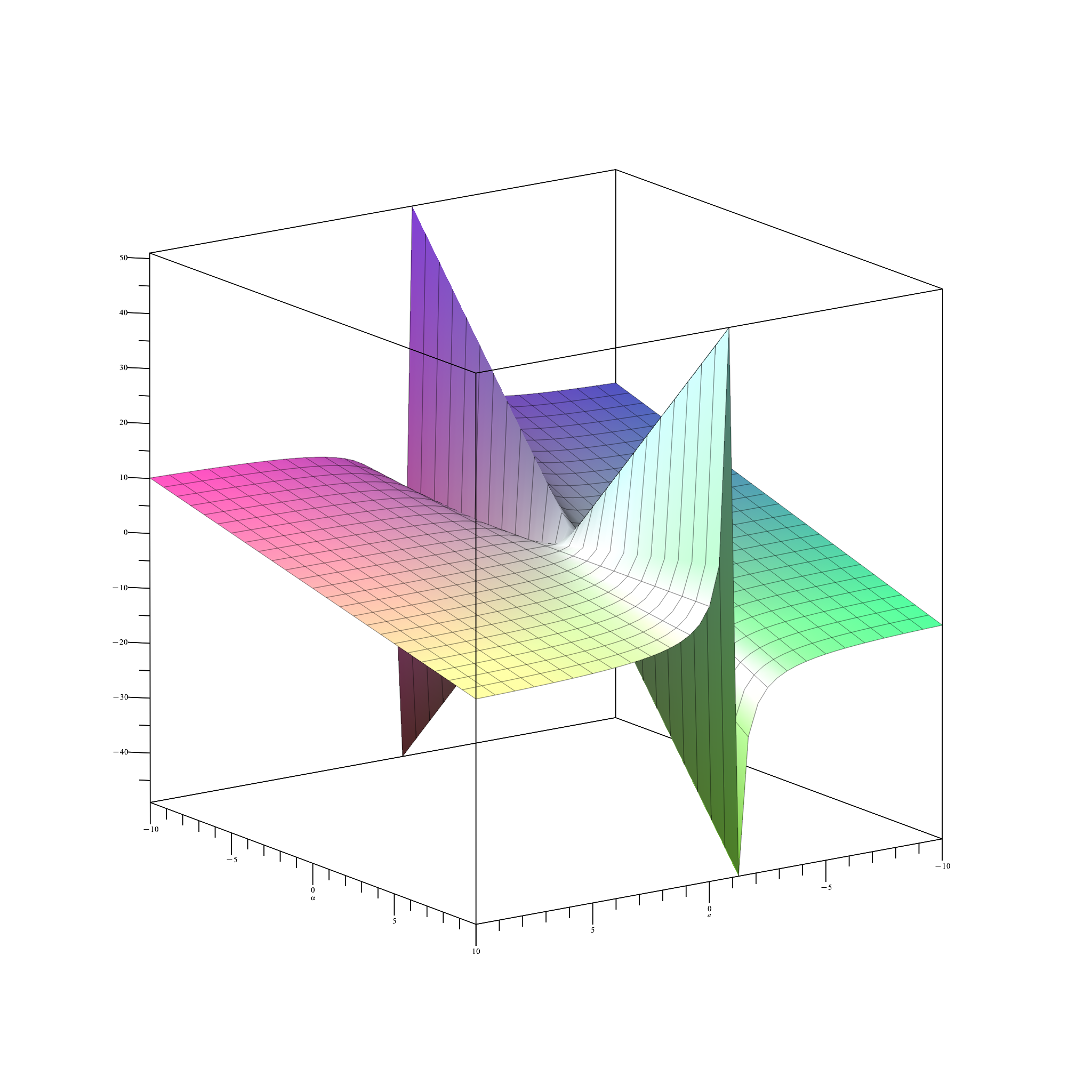}
       \vspace{-2cm}
       \begin{center}
       	\footnotesize{First Moment About the Origin}
       \end{center}
           \subsection{Second Moment about the Origin}

           \begin{equation*}
        \mu_{2}^{'} = E ( \xi^{2} )
        \end{equation*}
        \begin{equation*}
        =\int_{0}^{\infty} \xi^{2} F(\xi)d\xi
           \end{equation*}
\begin{equation*}
    =\int_{0}^{\infty}\xi^{2}(1+a)e^{-\xi}{\mathbb{E}}_{\alpha}(-a\xi^{\alpha})d\xi.
\end{equation*}
Using idempotent Representation
  \begin{equation*}
    = (1+a)\left[\left(\int_{0}^{\infty}\xi_{1}^{2}e^{-\xi_{1}}{\mathbb{E}}_{\alpha}(-a\xi_{1}^{\alpha})d\xi_{1}\right) e_{1} + \left(\int_{0}^{\infty}\xi_{2}^{2}e^{-\xi_{2}}{\mathbb{E}}_{\alpha}(-a\xi_{2}^{\alpha})d\xi_{2}\right)e_{2} \right]
     \end{equation*}
     \begin{equation*}
  = \sum_{k=0}^{\infty}(1+a)\left[\left( \int_{0}^{\infty}\xi_{1}^{2}e^{-\xi_{1}}\frac{(-a)^{k}\xi_{1}^{\alpha k} } {\Gamma(\alpha k+1)}d\xi_{1} \right)e_{1} + \left(\int_{0}^{\infty}\xi_{2}^{2}e^{-\xi_{2}}\frac{(-a)^{k}\xi_{2}^{\alpha k} }{\Gamma(\alpha k+1)}d\xi_{2}\right)e_{2}\right]
           \end{equation*}
            \begin{equation*}
     = \sum_{k=0}^{\infty}(1+a)\left[\left(\frac{(-a)^{k}}{\Gamma(\alpha k+1)} \Gamma(\alpha k+3)\right)e_{1}+ \left(\frac{(-a)^{k}}{\Gamma(\alpha k+1)} \Gamma(\alpha k+3)\right)e_{2}\right]
  \end{equation*}
  \begin{equation*}
 = (1+a)\sum_{k=0}^{\infty}[{(-a)^{k}(\alpha k+1) (\alpha k+2)} e_{1} + {(-a)^{k}(\alpha k+1)(\alpha k+2)} e_{2}]
  \end{equation*}
  \begin{equation*}
     = (1+a)\left[\left ( \frac{2}{1+a} - \frac{3 \alpha  a}{(1+a)^2 } +   \frac{ a \alpha^{2} ( a-1) }{(1+a)^{3}} \right) e_{1} + \left ( \frac{2}{1+a} - \frac{3 \alpha a}{(1+a)^2 } +  \frac{ a \alpha^{2} ( a-1) }{(1+a)^{3}}  \right) e_{2} \right]
  \end{equation*}
  \begin{equation*}
     = \left( 2 - \frac{3a \alpha }{1+a} + \frac{a(a-1) \alpha^{2}}{(1+a)^{2}}\right) (e_{1} +  e_{2} ).
  \end{equation*}
  
  Using (1.1) 
  \begin{equation}
   \mu_{2}^{'}  = 2 - \frac{3a \alpha}{1+a} + \frac{a(a-1) \alpha^{2}}{(1+a)^{2}}.
  \end{equation}
\mapleplot{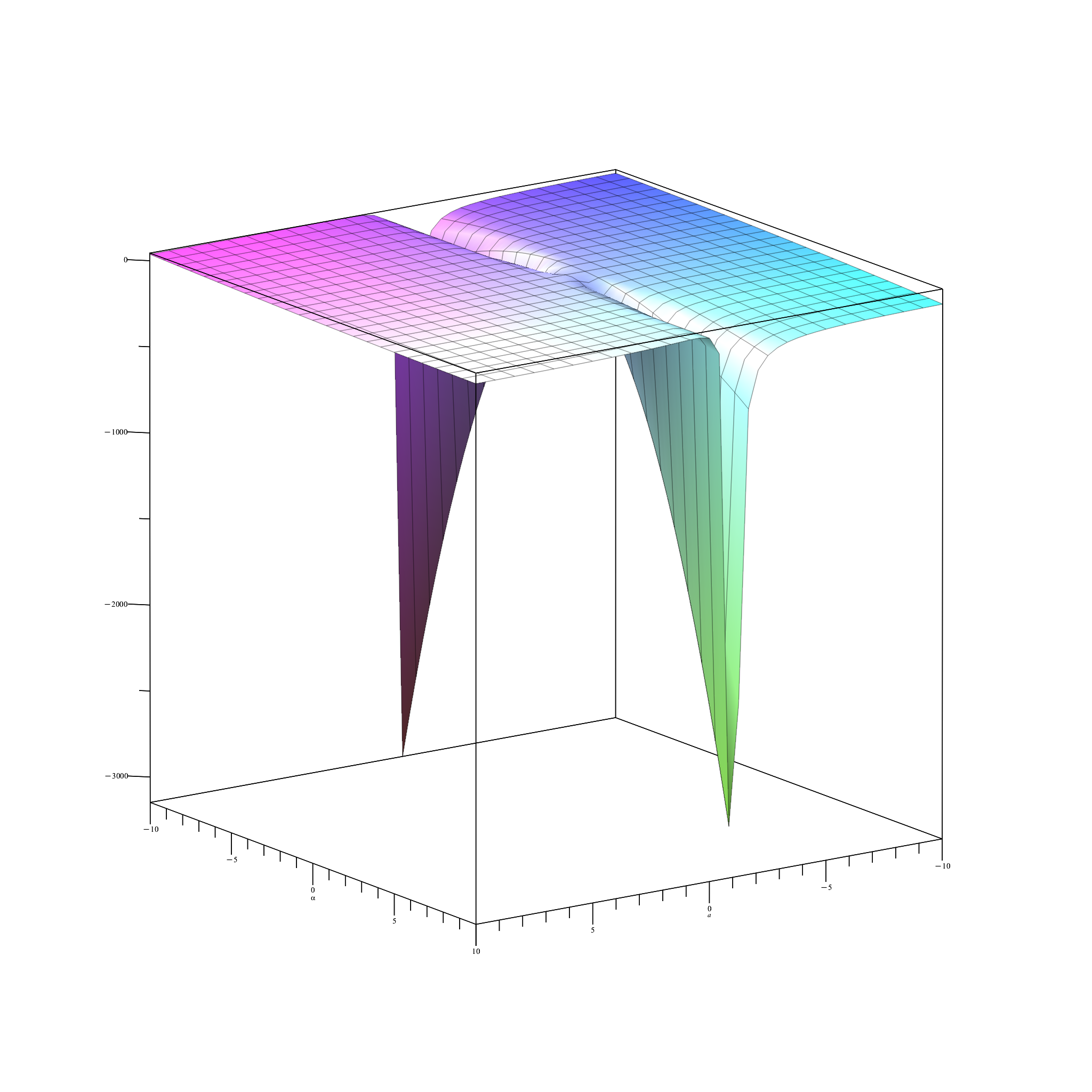}
 \vspace{-2cm}
\begin{center}
	\footnotesize{Second Moment About the Origin}
\end{center}
  \subsection{Third Moment about the Origin}

           \begin{equation*}
        \mu_{3}^{'} = E ( \xi^{3} )
        \end{equation*}
        \begin{equation*}
        =\int_{0}^{\infty} \xi^{3} F(\xi)d\xi
           \end{equation*}
\begin{equation*}
    =\int_{0}^{\infty}\xi^{3}(1+a)e^{-\xi}{\mathbb{E}}_{\alpha}(-a\xi^{\alpha})d\xi.
\end{equation*}
Using idempotent Representation
  \begin{equation*}
    = (1+a)\left[\left(\int_{0}^{\infty}\xi_{1}^{3}e^{-\xi_{1}}{\mathbb{E}}_{\alpha}(-a\xi_{1}^{\alpha})d\xi_{1}\right) e_{1} + \left(\int_{0}^{\infty}\xi_{2}^{3}e^{-\xi_{2}}{\mathbb{E}}_{\alpha}(-a\xi_{2}^{\alpha})d\xi_{2}\right)e_{2} \right]
     \end{equation*}
     \begin{equation*}
  = \sum_{k=0}^{\infty}(1+a)\left[\left( \int_{0}^{\infty}\xi_{1}^{3}e^{-\xi_{1}}\frac{(-a)^{k}\xi_{1}^{\alpha k} } {\Gamma(\alpha k+1)}d\xi_{1} \right)e_{1} + \left(\int_{0}^{\infty}\xi_{2}^{3}e^{-\xi_{2}}\frac{(-a)^{k}\xi_{2}^{\alpha k} }{\Gamma(\alpha k+1)}d\xi_{2}\right)e_{2}\right]
           \end{equation*}
            \begin{equation*}
     = \sum_{k=0}^{\infty}(1+a)\left[\left(\frac{(-a)^{k}}{\Gamma(\alpha k+1)} \Gamma(\alpha k+4)\right)e_{1}+ \left(\frac{(-a)^{k}}{\Gamma(\alpha k+1)} \Gamma(\alpha k+4)\right)e_{2}\right]
  \end{equation*}
  \begin{equation*}
 = (1+a)\sum_{k=0}^{\infty}[{(-a)^{k}(\alpha k+1) (\alpha k+2) (\alpha k+3)} e_{1} + {(-a)^{k}(\alpha k+1)(\alpha k+2) (\alpha k+3)} e_{2}]
  \end{equation*}
  \begin{align*}
     = (1+a) \left( \frac{6}{1+a} - \frac{11 a \alpha }{(1+a)^2 } + \frac{ 6a ( a-1)\alpha^{2}  }{(1+a)^{3}} - \frac{a(a^{2} -4a + 1)\alpha^{3}}{(1+a)^{4}} \right) e_{1} + \\ (1+a) \left( \frac{6}{1+a} - \frac{11 a \alpha }{(1+a)^2 } + \frac{ 6a ( a-1)\alpha^{2}  }{(1+a)^{3}} - \frac{a(a^{2} -4a + 1)\alpha^{3}}{(1+a)^{4}} \right) e_{2}
 \end{align*}
\begin{equation*}
      = \left( 6 - \frac{11 a \alpha}{(1+a) } + \frac{ 6a ( a-1)\alpha^{2} }{(1+a)^{2}} - \frac{a(a^{2} -4a + 1)\alpha^{3}}{(1+a)^{3}} \right)(e_{1} + e_{2}).
 \end{equation*}

 Using (1.1) 
 \begin{equation}
     \mu_{3}^{'} = 6 - \frac{11 a \alpha}{(1+a) } + \frac{ 6a ( a-1)\alpha^{2} }{(1+a)^{2}} - \frac{a(a^{2} -4a + 1)\alpha^{3}}{(1+a)^{3}}.
 \end{equation}
\mapleplot{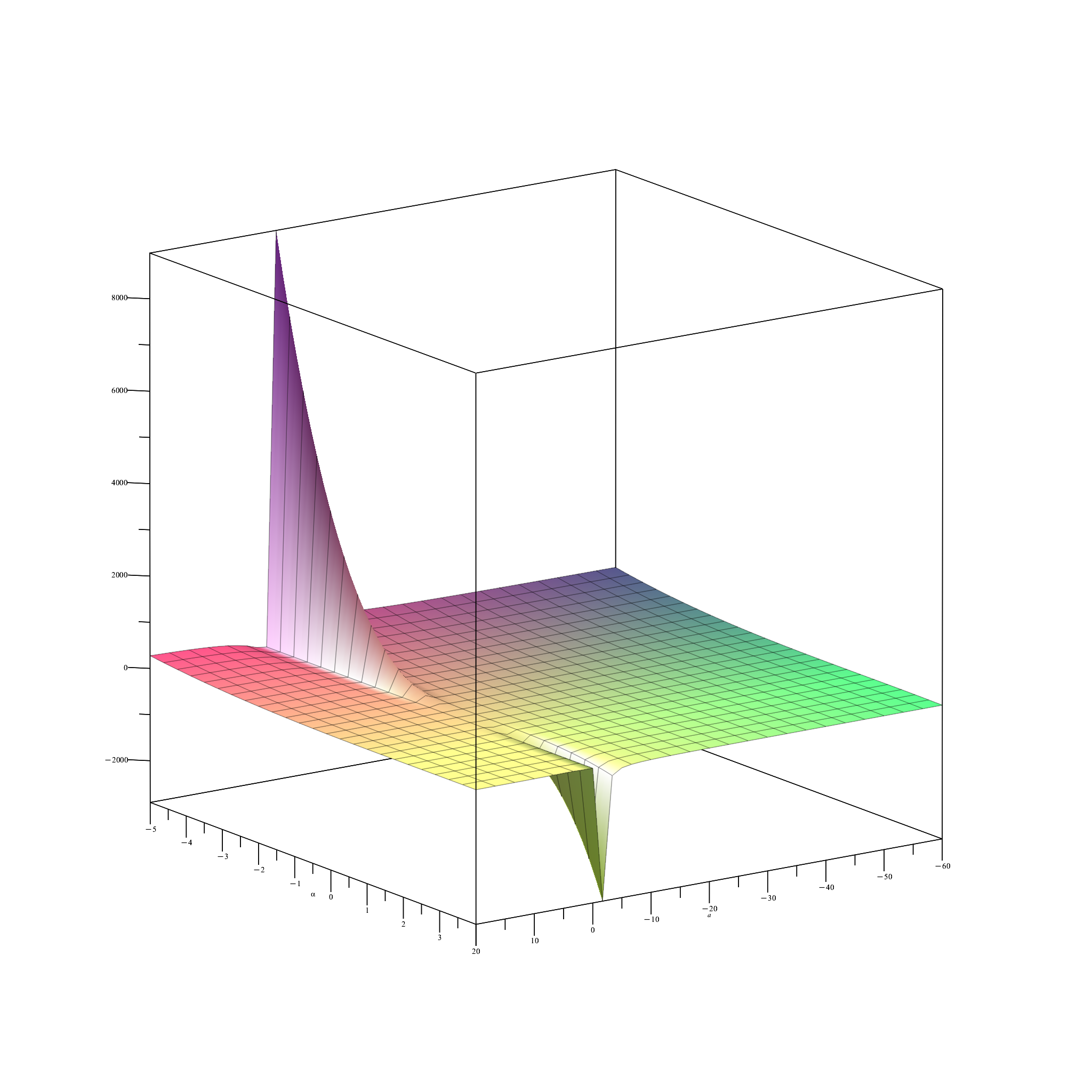}
 \vspace{-2cm}
\begin{center}
	\footnotesize{Third Moment About the Origin}
\end{center}
 \subsection{Fourth Moment about the Origin}

           \begin{equation*}
        \mu_{4}^{'} = E ( \xi^{4} )
        \end{equation*}
        \begin{equation*}
        =\int_{0}^{\infty} \xi^{4} F(\xi)d\xi
           \end{equation*}
\begin{equation*}
    =\int_{0}^{\infty}\xi^{4}(1+a)e^{-\xi}{\mathbb{E}}_{\alpha}(-a\xi^{\alpha})d\xi.
\end{equation*}
Using idempotent Representation
  \begin{equation*}
    = (1+a)\left[\left(\int_{0}^{\infty}\xi_{1}^{4}e^{-\xi_{1}}{\mathbb{E}}_{\alpha}(-a\xi_{1}^{\alpha})d\xi_{1}\right) e_{1} + \left(\int_{0}^{\infty}\xi_{2}^{4}e^{-\xi_{2}}{\mathbb{E}}_{\alpha}(-a\xi_{2}^{\alpha})d\xi_{2}\right)e_{2} \right]
     \end{equation*}
     \begin{equation*}
  = \sum_{k=0}^{\infty}(1+a)\left[\left( \int_{0}^{\infty}\xi_{1}^{4}e^{-\xi_{1}}\frac{(-a)^{k}\xi_{1}^{\alpha k} } {\Gamma(\alpha k+1)}d\xi_{1} \right)e_{1} + \left(\int_{0}^{\infty}\xi_{2}^{4}e^{-\xi_{2}}\frac{(-a)^{k}\xi_{2}^{\alpha k} }{\Gamma(\alpha k+1)}d\xi_{2}\right)e_{2}\right]
           \end{equation*}
            \begin{equation*}
     = \sum_{k=0}^{\infty}(1+a)\left[\left(\frac{(-a)^{k}}{\Gamma(\alpha k+1)} \Gamma(\alpha k+5)\right)e_{1}+ \left(\frac{(-a)^{k}}{\Gamma(\alpha k+1)} \Gamma(\alpha k+5)\right)e_{2}\right]
  \end{equation*}
  \begin{align*}
 = (1+a)\sum_{k=0}^{\infty}[{(-a)^{k}(\alpha k+1) (\alpha k+2) (\alpha k+3)(\alpha k+4)} e_{1} + \\ {(-a)^{k}(\alpha k+1)(\alpha k+2) (\alpha k+3)(\alpha k+4)} e_{2}]
  \end{align*}
 \begin{align*}
 =\left( 24-\frac{50a\alpha}{1+a}+\frac{35a(a-1)\alpha^{2}}{(1+a)^{2}}- \frac{10a(a^{2}-4a+1)\alpha^{3}}{(1+a)^{3}}+\frac{a(a^{3}-11a^{2}+11a-1)\alpha^{4}}{(1+a)^{4}}\right)e_{1}+\\ \left(24-\frac{50a\alpha}{1+a}+\frac{35a(a-1)\alpha^{2}}{(1+a)^{2}}-\frac{10a(a^{2}-4a+1)\alpha^{3}}{(1+a)^{3}}+\frac{a(a^{3}-11a^{2}+11a-1)\alpha^{4}}{(1+a)^{4}}\right)e_{2}
 \end{align*}
\begin{align*}
 =\left( 24-\frac{50a\alpha}{1+a}+\frac{35a(a-1)\alpha^{2}}{(1+a)^{2}}- \frac{10a(a^{2}-4a+1)\alpha^{3}}{(1+a)^{3}}+\frac{a(a^{3}-11a^{2}+11a-1)\alpha^{4}}{(1+a)^{4}}\right)(e_{1}+e_{2})
 \end{align*}
 Using (1.1)
\begin{align}
   \mu_{4}^{'}  = 24-\frac{50a\alpha}{1+a}+\frac{35a(a-1)\alpha^{2}}{(1+a)^{2}}- \frac{10a(a^{2}-4a+1)\alpha^{3}}{(1+a)^{3}}+\frac{a(a^{3}-11a^{2}+11a-1)\alpha^{4}}{(1+a)^{4}}.
\end{align}
\mapleplot{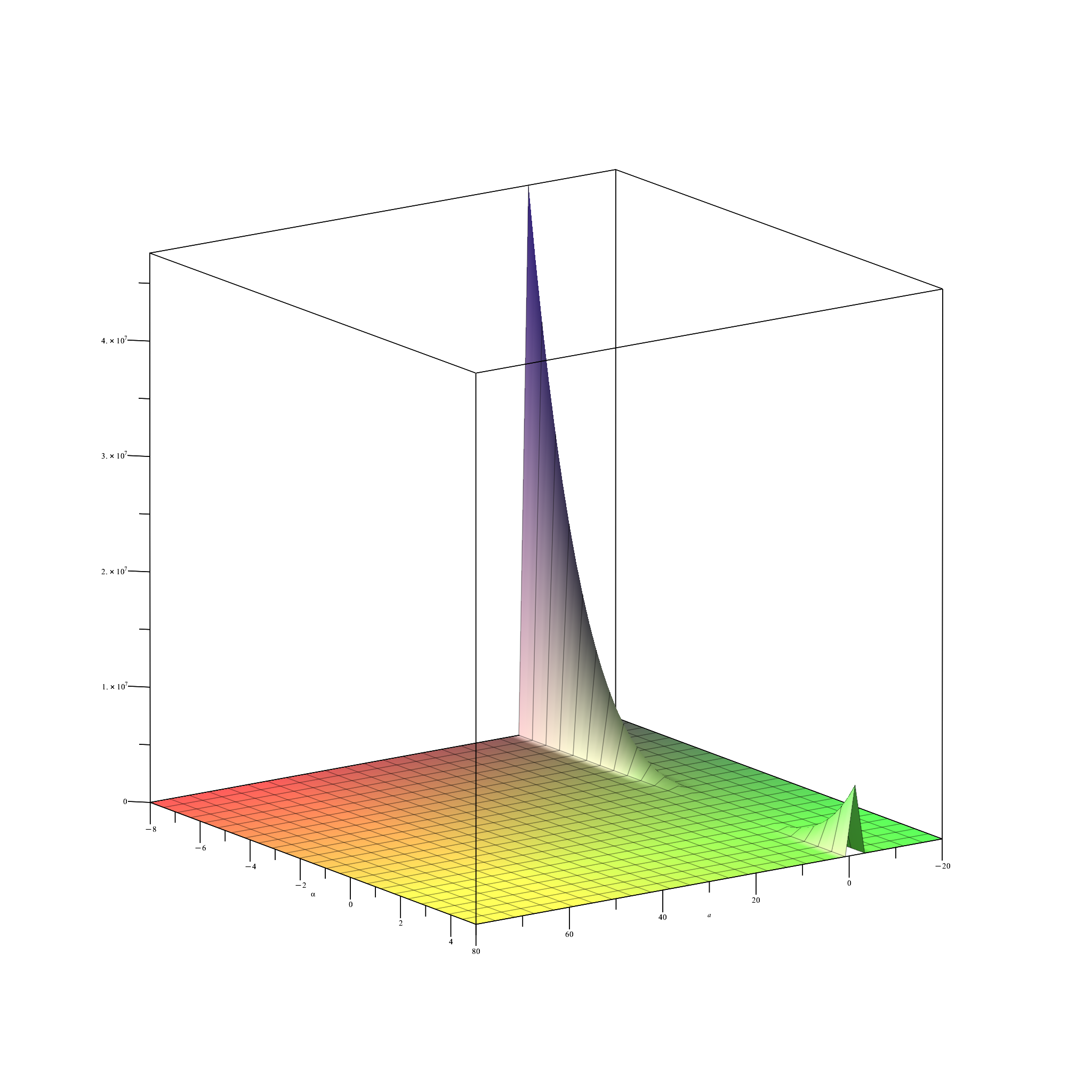}
 \vspace{-2cm}
\begin{center}
	\footnotesize{Fourth Moment About the Origin}
\end{center}
 \subsection{Mean and Variance}
   Mean of the probability density function (1.7)
      \begin{equation*}
     \text {Mean} =  1- \frac{a\alpha}{a+1}.
  \end{equation*}

\mapleinput

\mapleplot{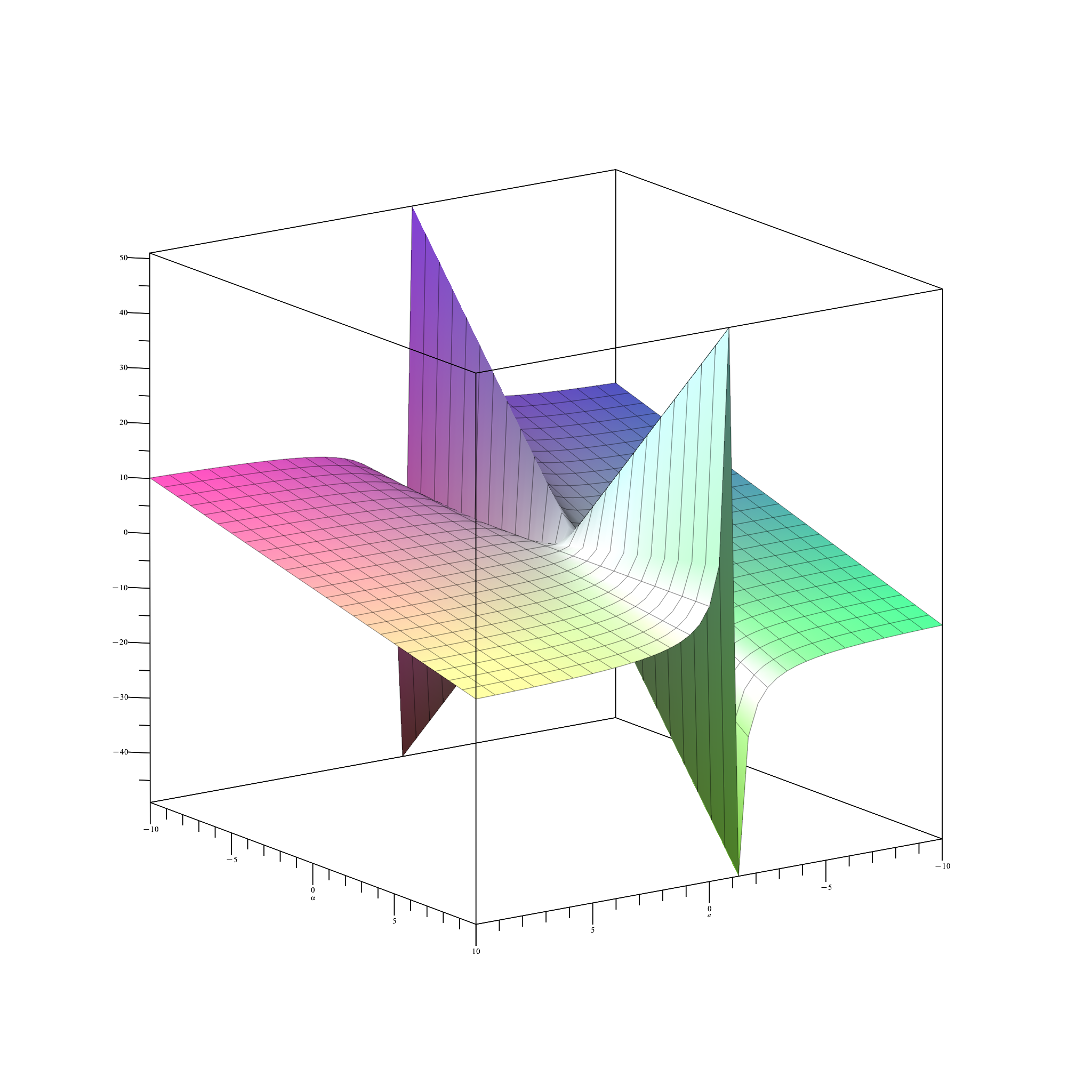}
 \vspace{-2cm}
\begin{center}
	\footnotesize{Mean}
\end{center}

The variance of the probability density function (1.7)
  \begin{equation*}
      \text{Variance} = 1 - \frac{a\alpha}{a+1} - \frac{a\alpha^{2}}{(a+1)^{2}}.
  \end{equation*}
\mapleplot{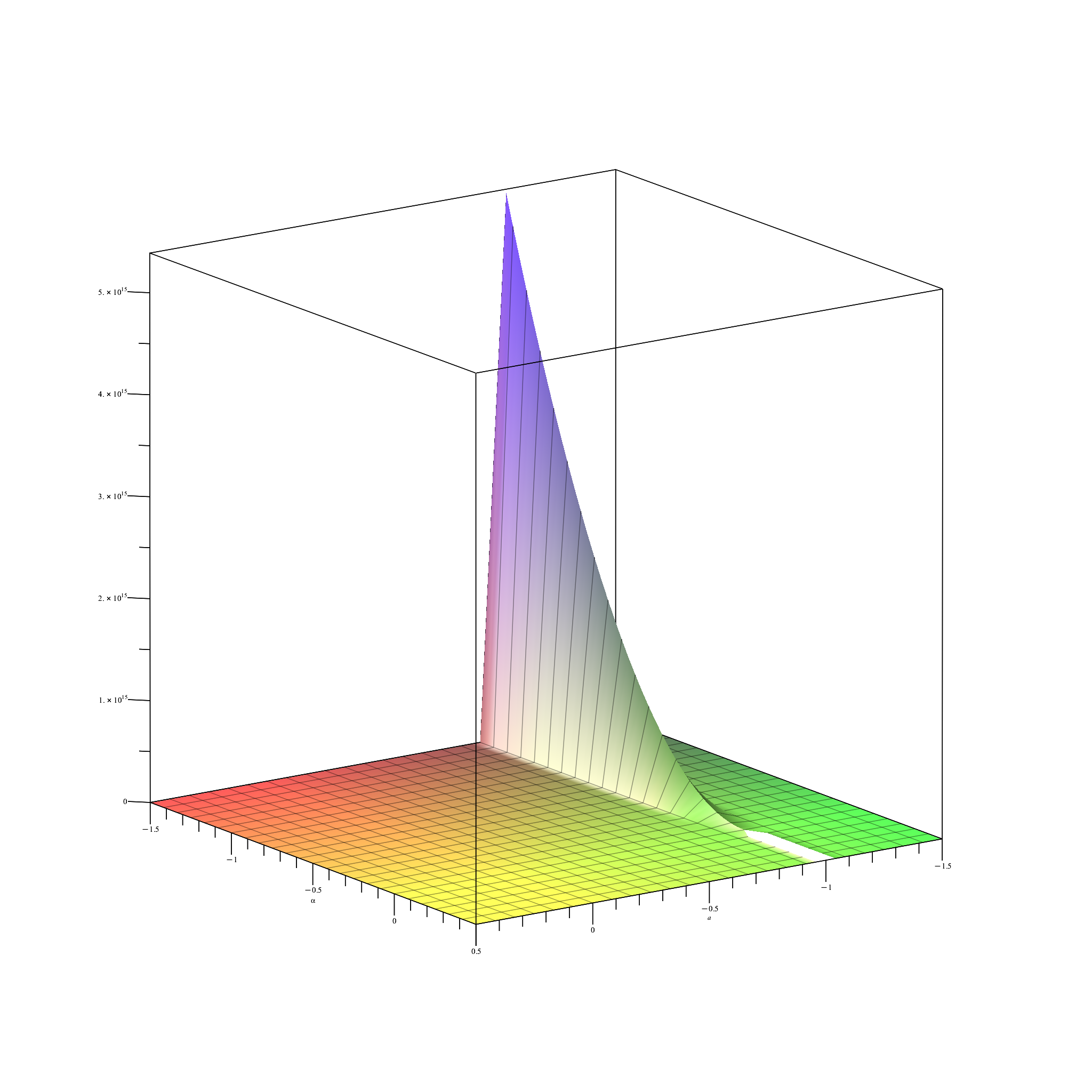}
 \vspace{-2cm}
\begin{center}
	\footnotesize{Variance}
\end{center}
{\bf{Conclusion.}} New results can be formed in the area of electromagnetism, quantum mechanics, and signal theory from the bicomplex Mittag-Leffler distribution and this will provide new avenues of research for researchers. \\ 
{\bf{Funding}} \\
There is no source of funding for this article\\
{\bf {Compliance with ethical standards}}\\
{\bf{Conflict of interest}} The author declares that this research involves no conflict of interest.\\
{\bf{Ethical approval}} This article does not contain any studies with human participants or animals.\\
{\bf {Availability of data and materials}}\\
 Not applicable

  \end{document}